\documentclass{amsart}
\usepackage{amsmath,amssymb}
\usepackage[UKenglish] {babel}

\usepackage{xy}
\xyoption{all}

  \def\R{\mathbb{R}} \def\C{\mathbb{C}}
\def\P{\mathbb{P}}

\def\O{\mathbb{O}}
\def\A{\mathbb{A}}
\def\H{\mathbb{H}}

\def\({\left(} \def\){\right)} 
\def\<{\langle} \def\>{\rangle}

\newcommand\forget[1]{}

\newenvironment{smatrix}{\left(\begin{smallmatrix}}{\end{smallmatrix}\right)}

\newtheorem{theorem}{Theorem}

\newtheorem{proposition}{Proposition}

\newtheorem{lemma}{Lemma}
\newtheorem{remarks}{Remarks}
\newtheorem{remarks and open problems}{Remarks and open problems}
\newtheorem{remark}{Remark}

\newtheorem{open problem}{Open problem}
\newtheorem{open problems}{Open problems}

\begin{document}

\large

\sloppy
\renewcommand{\thesection}{\arabic{section}}
\renewcommand{\theequation}{\thesection.\arabic{equation}}

\setcounter{page}{1}

\title{\bf Division algebras satisfying $(x^p, x^q, x^r)=0$}
\author{O. Diankha, A. Rochdi and M. Traor\'e}
\maketitle

\begin{abstract} We study algebras $A,$ over a field of characteristic zero, satisfying $(x^p, x^q, x^r)=0$ for $p, q, r$ in $\{1, 2\}.$ The existence of a unit element in such algebras leads to the third power-associativity. If, in addition, $A$ has degree $\leq 4$ then $A$ is power-commutative. We deduce that any $4$-dimensional real division algebra, with unit element, satisfying $(x^p, x^q, x^r)=0$ is quadratic. This persists for $(x, x^q, x^r)=0$ if we replace the word "unit" by "left-unit". \end{abstract}

\vspace{0.2cm} {\small Mathematics Subject Classification 2000: 17A05, 17A20, 17A30, 17A35, 17A45.}

\vspace{0.2cm}
{\small{\bf Keywords.} Third power-associative (Flexible, Quadratic) algebra. Division algebra.}

\vspace{0.6cm} \section{Introduction}

\vspace{0.3cm} The study of finite-dimensional (FD) real division algebras is a fascinating topic that arose after the discovery of quaternions $\H$ and octonions $\O.$ Classical results show that $\{\R, \C, \H\}$ classifies the FD real associative division algebras ([Fr 1878], [HKR 91]) and that $\{\R, \C, \H, \O\}$ classifies the FD real alternative\footnote{An algebra $A$ over an arbitrary field is said to be {\bf\em alternative} if it satisfies $(x, x, y)=(y, x, x)=0$ for all $x, y\in A.$} division algebras ([Zo 31], [HKR 91]).

\vspace{0.2cm} The $(1, 2, 4, 8)$-theorem ([Ho 40], [BM 58], [Ke 58], [HKR 91]), which states that $1,$ $2,$ $4,$ $8$ are the only possibilities for the dimension of every FD real division algebra, revolutionized the theory and has become an indispensable tool in almost all subsequent work. On the other hand the classification of real division algebras is trivial in dimension $1$ ([S 66] p. 2), ([Rod 04] Proposition 1.2 p. 104). It was completed and refined quite recently in dimension $2$ ([HP 04], [Di 05]) after a few attempts ([BBO 81, 82], [AK 83], [Bu 85], [G 98]). Thus work focused onto algebras of dimension $4$ and $8.$

\vspace{0.2cm} Studies in [BBO 82], where the pseudo-octonions algebra $\P$ [Ok 78] played an important role, led to the classification of the FD real flexible\footnote{Algebra $A$ is said to be {\bf\em flexible} if it satisfies $(x, y, x)=0$ for all $x, y\in A.$} division algebras ([CDKR 99], [Da 06]) generalizing [Zo 31].

\vspace{0.2cm} Another generalizing result ([Os 62], [Di 00]) classifies the $4$-dimensional real quadratic\footnote{Algebra $A$ is said to be {\bf\em quadratic} if it contains a unit element $e$ and for all $x\in A$ the elements $e,$ $x,$ $x^2$ are linearly dependant.} division algebras, and despite a series of additional works in dimension $8$ ([DL 03], [Li 04], [DFL 06], [DR 10]) the classification of $8$-dimensional real quadratic division algebras remains an open problem.

\vspace{0.2cm} The class of power-commutative algebras contains the flexible [Raf $50_2$], and quadratic algebras. Note here that among the finite-dimensional real division algebras, power-associative algebras are the same as the quadratic algebras ([Roc 94] Corollaire 2.27 p. 36). The problem of determining all real division algebras is more difficult in dimension $8$ than on dimension $4$ where the works are progressing relatively more quickly, as occurred for the power-commutative algebras [DR 11].

\vspace{0.2cm} Third power-associative algebras are a natural generalization of power-commutative ones, and was already studied ([Elm 83, 87], [EP 94]) but apparently never in the context of division algebras. In the present work we begin our study with a third power-associative algebra over an arbitrary field of characteristic zero. We show that the additional assumption of the existence of a one-sided unit and that the algebra has degree $\leq 4,$ ensures the commutativity of powers (Lemma 1). We deduce that every real third power-associative division algebra, with one-sided unit, having degree $\leq 4$ is quadratic (Theorem 3).

\vspace{0.2cm} Algebras satisfying the identity $(x^p, x^q, x^r)=0$ for fixed $p, q, r\in\{1, 2\}$ were also studied outside the general context of division algebras ([CR 08], [CRR 11], [Elm 01], [EE 04]). Here we show first that every real division algebra, with left-unit, satisfying $(x,x,x^2)=0$ is quadratic (Theorem 1). This remains true for every real division algebra, with unit element, of degree $\leq 4$ satisfying $(x^p, x^q, x^r)=0$ (Theorem 2). Above result persists for every identity $(x, x^q, x^r)=0$ if we replace the word "unit element" by "left-unit element" (Theorem 3). In the other hand there are examples of real left-unital division algebras of degree $2$ satisfying $(x^2, x^q, x^r)=0$ for all $q, r\in\{1, 2\}$ containing no right-unit element. Hence the need, here, for the assumption of the existence of a two-sided unit element.

\vspace{0.7cm} \section{Notations and preliminary results}

\vspace{0.3cm} Algebras $A$ will be considered over a field $\mathbb{K}$ of characteristic zero. We denote by $(x, y, z)$ (resp. $[x, y]$) the associator $(xy)z-x(yz)$ of $x, y, z\in A$ (resp. the commutator of $x, y\in A$). The subalgebra of $A$ generated by every element $x\in A$ is denoted by $A(x).$

\vspace{0.2cm} \begin{enumerate} \item Algebra $A$ is called

\vspace{0.2cm} \begin{enumerate} \item {\bf\em third power-associative} {\bf (TPA)} if $(x, x, x)=0$ for all $x\in A,$

\vspace{0.1cm} \item {\bf\em power-associative} if $A(x)$ is associative for all $x\in A,$

\vspace{0.1cm} \item {\bf\em power-commutative} if $A(x)$ is commutative for all $x\in A.$ \end{enumerate}

\vspace{0.2cm} \item If $A$ is FD then the set $\{dim(A(y)): y\in A\}$ admit a bigger finite element $d$ called the degree of $A$ [Rod 94].

\vspace{0.2cm} \item $A$ is called a division algebra if it is FD and the operators $L_x$ and $R_x$ of left and right multiplication by $x$ are bijective for all $x\in A-\{0\}.$

\vspace{0.2cm} \item If $\mathbb{K}$ is the field of real numbers $\R$ the algebra is said to be real. The $(1, 2, 4, 8)$-theorem shows that the degree of every real division algebra is $1, 2, 4$ or $8.$ $\Box$ \end{enumerate}

\vspace{0.3cm} Let now $\A$ be one of classical real division algebras $\C$ (complex numbers), $\H$ (quaternions) or $\O$ (octonions), and $^*\A,$ $\stackrel{*}{\A},$ the standard isotopes of $\A$ having $\A$ as vectorial space and products $x*y$ given respectively by $\overline{x}y,$ $\overline{x}\ \overline{y}$ where $x\mapsto\overline{x}$ is the standard conjugation of $\A.$

\vspace{0.7cm} \section{The identities $(x^p, x^q, x^r)=0$}

\vspace{0.3cm} Let now $A$ be an algebra and $p, q, r$ be natural numbers in $\{1, 2\}.$ There are maps $f_m:A\rightarrow A$ for $m=1, \dots, p+q+r-1$ such that the equality

\vspace{0.1cm} \begin{eqnarray*} \Big( (x+\lambda y)^p, (x+\lambda y)^q, (x+\lambda y)^r\Big) &=& (x^p, x^q, x^r)+\lambda f_1(x, y)+\dots \\ && \quad +\lambda^{p+q+r-1}f_{p+q+r-1}(x, y) +\lambda^{p+q+r}(y^p, y^q, y^r) \end{eqnarray*}

\vspace{0.1cm} holds for all scalar $\lambda$ and $x, y\in A.$ Moreover, $f_m(y,x)=f_{p+q+r-m}(x,y)$ for all $m\in\{1, \dots, p+q+r-1\}.$ The identity $(x^p, x^q, x^r)=0$ in $A$ is equivalent to $f_m\equiv 0$ for all $m.$ We denote $f_m\equiv 0$ by (p.q.r.m) and call it the $m^{th}$ identity obtained by linearization of $(x^p, x^q, x^r)=0.$ The following three tables, where $x\bullet y$ denotes $xy+yx,$ specify the cases (p.q.r.1), (p.q.r.2) and (2.2.2.3) respectively:

\[ \begin{tabular}{cc} \\
\hline \multicolumn{1}{|c|}{$(x^p,x^q,x^r)=0$} &
\multicolumn{1}{|c|}{The first corresponding identity (p.q.r.1)} \\ \hline
\multicolumn{1}{|c|}{$(x,x,x)=0$} & \multicolumn{1}{|c|}{$[x^2,y]+[x\bullet y,x]=0$} \\ \hline
\multicolumn{1}{|c|}{$(x,x,x^2)=0$} & \multicolumn{1}{|c|}{$(x,x,x\bullet y)+(x,y,x^2)+(y,x,x^2)=0$} \\ \hline
\multicolumn{1}{|c|}{$(x,x^2,x)=0$} & \multicolumn{1}{|c|}{$(x,x^2,y)+(x,x\bullet y,x)+(y,x^2,x)=0$} \\ \hline
\multicolumn{1}{|c|}{$(x,x^2,x^2)=0$} & \multicolumn{1}{|c|}{$(x,x^2,x\bullet y)+(x,x\bullet y,x^2)+(y,x^2,x^2)=0$} \\ \hline \multicolumn{1}{|c|}{$(x^2,x,x)=0$} & \multicolumn{1}{|c|}{$(x^2,x,y)+(x^2,y,x)+(x\bullet y,x,x)=0$} \\ \hline
\multicolumn{1}{|c|}{$(x^2,x,x^2)=0$} & \multicolumn{1}{|c|}{$(x^2,x,x\bullet y)+(x^2,y,x^2)+(x\bullet y,x,x^2)=0$} \\ \hline \multicolumn{1}{|c|}{$(x^2,x^2,x)=0$} & \multicolumn{1}{|c|}{$(x^2,x^2,y)+(x^2,x\bullet y,x)+(x\bullet y,x^2,x)=0$} \\ \hline \multicolumn{1}{|c|}{$(x^2,x^2,x^2)=0$} & \multicolumn{1}{|c|}{$(x^2,x^2,x\bullet y)+(x^2,x\bullet y,x^2)+(x\bullet y,x^2,x^2)=0$} \\ \hline \end{tabular} \]
\[ \begin{tabular}{cc} \\
\hline \multicolumn{1}{|c|}{$(x^p,x^q,x^r)=0$} &
\multicolumn{1}{|c|}{The second corresponding identity (p.q.r.2)} \\ \hline
\multicolumn{1}{|c|}{$(x,x,x)=0$} & \multicolumn{1}{|c|}{$[x^2,y]+[x\bullet y,x]=0$} \\ \hline
\multicolumn{1}{|c|}{$(x,x,x^2)=0$} & \multicolumn{1}{|c|}{$(x,x,y^2)+(x,y,x\bullet y)+(y,x,x\bullet y)+(y,y,x^2)$} \\ \hline
\multicolumn{1}{|c|}{$(x,x^2,x)=0$} & \multicolumn{1}{|c|}{$(x,x\bullet y,y)+(x,y^2,x)+(y,x^2,y)+(y,x\bullet y,x)=0$} \\ \hline
\multicolumn{1}{|c|}{$(x,x^2,x^2)=0$} & \multicolumn{1}{|c|}{$(x,x^2,y^2)+(x,y^2,x^2)+(y,x\bullet y,x^2)+(y,x^2,x\bullet y)$} \\
\multicolumn{1}{|c|}{} & \multicolumn{1}{|c|}{$+(x,x\bullet y,x\bullet y)=0$} \\ \hline
\multicolumn{1}{|c|}{$(x^2,x,x)=0$} & \multicolumn{1}{|c|}{$(x^2,y,y)+(x\bullet y,x,y)+(x\bullet y,y,x)+(y^2,x,x)=0$} \\ \hline
\multicolumn{1}{|c|}{$(x^2,x,x^2)=0$} & \multicolumn{1}{|c|}{$(x^2,x,y^2)+(x^2,y,x\bullet y)+(x\bullet y,x,x\bullet y)$} \\
\multicolumn{1}{|c|}{} & \multicolumn{1}{|c|}{$+(x\bullet y,y,x^2)+(y^2,x,x^2)=0$} \\ \hline
\multicolumn{1}{|c|}{$(x^2,x^2,x)=0$} & \multicolumn{1}{|c|}{$(x^2,x\bullet y,y)+(x^2,y^2,x)+(x\bullet y,x^2,y)$} \\
\multicolumn{1}{|c|}{} & \multicolumn{1}{|c|}{$+(x\bullet y,x\bullet y,x)+(y^2,x^2,x)=0$} \\ \hline
\multicolumn{1}{|c|}{$(x^2,x^2,x^2)=0$} &  \multicolumn{1}{|c|}{$(x^2,x^2,y^2)+(x^2,x\bullet y,x\bullet y)+(x^2,y^2,x^2)$} \\ \multicolumn{1}{|c|}{} & \multicolumn{1}{|c|}{$+(x\bullet y,x^2,x\bullet y)+(x\bullet y,x\bullet y,x^2)+(y^2,x^2,x^2)=0$} \\ \hline \end{tabular} \]
\[ \begin{tabular}{cc} \\ \hline
\multicolumn{1}{|c|}{} & \multicolumn{1}{|c|}{The third corresponding identity (2.2.2.3)} \\ \cline{2-2} \multicolumn{1}{|c|}{$(x^2,x^2,x^2)=0$} & \multicolumn{1}{|c|}{$(x\bullet y,x^2,y^2)+(x^2,x\bullet y,y^2)+(x^2,y^2,x\bullet y)+(x\bullet y,y^2,x^2)$} \\ \multicolumn{1}{|c|}{} & \multicolumn{1}{|c|}{$+(y^2,x\bullet y,x^2)+(y^2,x^2,x\bullet y)+(x\bullet y,x\bullet y,x\bullet y)$} \\ \hline
\end{tabular} \]

\vspace{0.5cm} We have some relationships between the identities $(x^p, x^q, x^r)=0:$

\vspace{0.3cm} \begin{proposition} Every TPA algebra satisfies the identity $(x, x^2, x)=0.$  \end{proposition}

\vspace{0.2cm} {\bf Proof.} The result is well known [Raf $50_1$] and follows here from (1.1.1.1) by putting $y=x^2$ and taking into account the characteristic of $\mathbb{K}.$ $\Box$

\vspace{0.3cm} \begin{remark} Third power-associativity is stronger than identity $(x, x^2, x)=0$ {\em ([Ch 09] Remarque 1.17, p. 18)}. However, we know no examples of division algebras satisfying the identity $(x, x^2, x)=0$ which are not TPA. $\Box$ \end{remark}

\vspace{0.3cm} \begin{proposition} Let $A$ be an algebra, with unit element $e,$ satisfying the identity $(x^p, x^q, x^r)=0$ for fixed $p, q, r$ in $\{1, 2\}.$ Then $A$ is TPA.  \end{proposition}

\vspace{0.2cm} {\bf Proof.} We can assume that $(p, q, r)\neq (1, 1, 1)$ and the result is immediately obtained from the equality (p.q.r.m), where $m=p+q+r-3,$ by putting $y=e$ and taking into account the characteristic of $\mathbb{K}.$ $\Box$

\vspace{0.5cm} \begin{proposition} Let $A$ be an algebra over $\mathbb{K},$ with a left-unit $e,$ having no non-zero divisors of zero. If $A$ satisfy the identity $(x, x^q, x^r)=0$ for fixed $q, r$ in $\{1, 2\},$ then $A$ has unit element $e$ and is TPA. \end{proposition}

\vspace{0.2cm} {\bf Proof.} We can assume that $(q, r)\neq (1, 1)$ and putting $x=e$ in the equality (1.q.r.1), we have: \ $0=(y, e, e)=(ye-y)e.$ As $A$ has no non-zero divisors of zero, $e$ is a right-unit and Proposition 2 conclude. $\Box$

\vspace{0.3cm} The identities $(x, x, x^2)=0,$ $(x^2, x, x)=0$ have a particular importance:

\vspace{0.3cm} \begin{proposition} Let $A$ be an algebra over $\mathbb{K},$ with left-unit $e,$ having no non-zero divisors of zero. If $A$ satisfy the identity $(x, x, x^2)=0$ then $A$ has unit element $e$ and is power-associative. \end{proposition}

\vspace{0.2cm} {\bf Proof.} Consequence of Proposition 3 and ([A 48] Theorem 2). $\Box$

\vspace{0.3cm} As consequence:

\vspace{0.3cm} \begin{theorem} Let $A$ be a real division algebra, with left-unit $e,$ satisfying the identity $(x, x, x^2)=0.$ Then $A$ has unit element $e$ and is a quadratic algebra. \end{theorem}

\vspace{0.2cm} {\bf Proof.} $A$ is power-associative and [Roc 94] conclude. $\Box$

\vspace{0.3cm} We illustrate by the organization chart summarized below the "hierarchy" between different non-associative identities introduced in this text:

\vspace{0.2cm} {\Large\bf \begin{eqnarray*} \mbox{ Associativity } \hspace{2.8cm} . \\  \hspace{3.5cm} \Downarrow\hspace{4cm} . \\
\mbox{ Alternativity } \hspace{2.8cm} . \\
\Downarrow\hspace{1.5cm}\Downarrow\hspace{3cm} . \\
\mbox{ Flexiblility } \hspace{0.3cm} \mbox{ Power-associativity } \\
\Downarrow\hspace{1.5cm}\Downarrow\hspace{3cm} . \\
\mbox{ Power-commutativity } \hspace{1.7cm} . \\
\Downarrow \hspace{4cm} . \\
\mbox{ Third power-associativity } \hspace{1.3cm} . \\
\Downarrow \hspace{4cm} . \\
\mbox{ $(x, x^2, x)=(x^2, x^2, x^2)=0$ } \hspace{1.6cm} . \end{eqnarray*}}

\vspace{0.8cm} \section{The result}

\vspace{0.5cm} We have the following preliminary result:

\vspace{0.3cm} \begin{lemma} Let $A$ be an algebra, over $\mathbb{K},$ of degree $\leq 4$ having a left unit $e.$ Then $A$ is power-commutative in any one of the following cases:

\vspace{0.3cm} \begin{enumerate} \item $e$ is a two-sided unit and $A$ satisfies the identity $(x^p, x^q, x^r)=0$ for fixed $p, q, r$ in $\{1, 2\}.$

\vspace{0.2cm} \item $A$ has no non-zero divisor of zero and satisfies the identity $(x, x^q, x^r)=0$ for fixed $q, r$ in $\{1, 2\}.$ \end{enumerate} \end{lemma}

\vspace{0.3cm} {\bf Proof.} {\bf (1)} According to the Proposition 2 algebra $A$ is TPA. The first assertion follows then from ([Raf 50] Th\'eor\`eme 3, p. 578).

\vspace{0.2cm} {\bf (2)} The result follows from the first assertion and Proposition 3. $\Box$

\vspace{0.5cm} We can now state the following result:

\vspace{0.3cm} \begin{theorem} Let $A$ be a real division algebra, with unit element $e,$ of degree $\leq 4.$ Then the following assertions are equivalent:

\vspace{0.2cm} \begin{enumerate} \item $A$ satisfy the identity $(x^p, x^q, x^r)=0$ for fixed $p, q, r$ in $\{1, 2\}.$

\vspace{0.2cm} \item $A$ is power-associative.

\vspace{0.2cm} \item $A$ is quadratic. \end{enumerate} \end{theorem}

\vspace{0.3cm} {\bf Proof.} {\bf (1) $\Rightarrow$ (2).} $A$ is power-commutative by Lemma 1. The result is then a consequence of Hopf's commutative theorem [H 40] and Yang-Petro's theorem ([Y 81], [Pet 87]).

\vspace{0.2cm} {\bf (2) $\Rightarrow$ (3).} See ([Roc 94] Corollaire 2.27 p. 36).

\vspace{0.2cm} {\bf (3) $\Rightarrow$ (1)} is obvious. $\Box$

\vspace{0.5cm} Theorem 2 remain true if one replaces the assumption "$A$ has a unit element" by "$A$ has a left-unit element" and identity "$(x^p, x^q, x^r)=0$" by "$(x, x^q, x^r)=0$", thanks to Proposition 3:

\vspace{0.3cm} \begin{theorem} Let $A$ be a real division algebra of degree $\leq 4$ having a left-unit element. Then the following assertions are equivalent:

\vspace{0.2cm} \begin{enumerate} \item $A$ satisfy the identity $(x, x^q, x^r)=0$ for fixed $q, r$ in $\{1, 2\}.$

\vspace{0.2cm} \item $A$ is quadratic. $\Box$ \end{enumerate} \end{theorem}

\vspace{0.3cm} \begin{remarks} Let $\A$ any one of real division algebras $\C,$ $\H,$ $\O.$

\vspace{0.2cm} \begin{enumerate} \item The isotope standard $\stackrel{*}{\A}$ of $\A$ and the pseudo-octonion algebra $\P$ are flexible {\bf\em but not power-associative} algebras. Thus the hypothesis of the existence of a left-unit element in Theorem 3 is necessary.

\vspace{0.2cm} \item The division algebra $^*{\A},$ having left-unit and degree $2,$ satisfies all the identities $(x^2, x^q, x^r)=0$ but is {\bf\em not TPA}. So Proposition 3 and Theorem 3 does not have a similar for the identities: $(x^2, x^q, x^r)=0.$

\vspace{0.2cm} \item Taking into account Theorem 2 and the $(1, 2, 4, 8)$-theorem it is natural to wonder if there are examples of unital TPA real division algebras of degree $8.$ $\Box$ \end{enumerate} \end{remarks}

\vspace{0.5cm}
{\small O. Diankha

\vspace{0.1cm} D\'epartement de Math\'ematiques et Informatique, Facult\'e des Sciences et Techniques, Universit\'e Cheikh Anta Diop, Dakar, Senegal

\vspace{0.1cm} e-mail: odiankha@ucad.sn

\vspace{0.3cm} A. Rochdi and M. Traor\'e

\vspace{0.1cm} D\'epartement de Math\'ematiques et Informatique, Facult\'e des Sciences Ben M'Sik, Universit\'e Hassan II-Mohammedia, 7955 Casablanca, Morocco

\vspace{0.1cm} e-mail: abdellatifro@yahoo.fr and sasmohasas@yahoo.fr}

\end{document}